\documentclass[12pt]{article}
\usepackage{amsmath,amssymb,amssymb}

\newtheorem{lemma}{Lemma}
\newtheorem{theorem}{Theorem}
\newtheorem{definition}{Definition}
\newenvironment{proof}{
  \noindent\textbf{Proof}\ }{\hspace*{\fill}
  \begin{math}\Box\end{math}\medskip}
\newenvironment{proof theorem}{
  \noindent\textbf{Proof of Theorem 1}\ }{\hspace*{\fill}
  \begin{math}\Box\end{math}\medskip}

\begin{document}

\title{Two identities of derangements} \author{Le Anh Vinh\\
School of Mathematics\\University of New South Wales\\Sydney 2052 NSW} 
\date{\empty}
\maketitle

\begin{abstract}
 In this note, we present two new identities for derangements. As a corrolary, we have a combinatorial proof of the irreducibility of the standard representation of symmetric groups.
\end{abstract}

\section{Introduction}
A derangement is the permutation $\sigma$ of $\{1,2,\ldots,n\}$ that there is no $i$ satisfying $\sigma(i) = i$. It is well-known that the number $d(n)$ of dereangements equals:
\[ d(n) = \sum_{k=0}^n (-1)^k \frac{n!}{k!}\]
and satisfies the following identity (since both sides are the number of permutations on $n$ letters)
\begin{equation}\label{easy}\sum_{k = 0}^n \binom{n}{k} d ( k ) = n!.\end{equation}

The Stirling set number $S(n,m)$ is the number of ways of partitioning a set of $n$ elements into $m$ nonempty sets. We definte $[x]_r = x(x-1)\ldots(x-r+1)$ (by convention $[x]_0=1$). Then (see \cite{stirling math})
\begin{equation}\label{stirling} x^n = \sum_{m=0}^n S(n,m)[x]_r.\end{equation}

The number of ways a set of  elements can be partitioned into nonempty subsets is called a Bell number and is denoted $B_n$. We use the convention that $B_0 = 1$. The integer $B_n$ can also be define by the sum (see \cite{bell math})
\begin{equation}\label{bell} B_n = \sum_{m=0}^n S(n,m)\end{equation}

The main results of this note are the following generalizations of (\ref{easy}). 

\begin{theorem}\label{main} Let $n, k, l$ be three natural numbers. Then
\begin{equation} \label{hard}
\sum_{k = 0}^n \binom{n-l}{k-l} d(n-k) = ( n - l ) !. 
\end{equation}
\end{theorem}

\begin{theorem}
Suppose that $n \geq m$ are two natural numbers. Let $g(x) = a_m x^m +\ldots+ a_0$ be a polynomial with integer coefficients. Then
\begin{equation}\label{second}\sum_k g(k)\binom{n}{k} d(n-k) = \left\{\sum_{i=0}^m a_i B_i \right\} n!.\end{equation}
\end{theorem}

We use the convention that $\binom{n}{m} = 0$ if $m < 0$ or $n < m$. Also set $d( k ) = 0$ if $k < 0$  and $d(0) = 1$. Note that taking $l = 0$ in (\ref{hard}) implies (\ref{easy}) since $\binom{n}{k} = \binom{n}{n-k}$.

\section{Some Lemmas}

We define $f_n( k )$ to be the number of permutations of $\{ 1, \ldots, n \}$ that fix exactly $k$
positions. By convention, $f_n ( k ) = 0$ if $k < 0$ or $k >n$. We have the  following recursion for $f_n ( k )$.

\begin{lemma}\label{re} Suppose that $n, k$ are positive integers. Then
  \[f_{n + 1} ( k ) = f_n ( k - 1 ) + ( n - k ) f_n ( k ) + ( k + 1 ) f_n ( k +1 ).\] 
\end{lemma} 

\begin{proof}
  Let $\sigma$ be any permutation of $\{ 1, \ldots, n + 1 \}$ which has
  exactly $k$ fixed points. We have two cases.
  
\begin{enumerate}
	\item Suppose that $\sigma ( n + 1 ) = n + 1$. Then $\sigma$ corresponds
  to a restricted permutation on $\{ 1, \ldots, n\}$ which fixes $k - 1$
  points of $\{ 1, \ldots, n\}$.
  \item Suppose that $\sigma ( n + 1 ) = i$ for some $i \in \{ 1, \ldots, n\}$. Then there exists 
  $j \in \{ 1, \ldots, n \}$ such that $\sigma ( j ) = n+ 1$. There are two separate subcases.
  
\begin{enumerate}
	\item If $i = j$ then we can obtain a correspondence between $\sigma$ and a 
  permutation $\sigma'$ of $\{ 1, \ldots, n \}$ from $\sigma$ as follows: $\sigma' ( i ) = i$ and $\sigma' (
  t ) = \sigma ( t )$ for $t \neq i$. It is clear that $\sigma'$ has $k + 1$ fixed points.
  Conversely, for each permutation of $\{ 1, \ldots, n \}$ that has $k + 1$
  fixed points, we can choose $i$ to be any of its fixed points and then swapping $i$ and $n+1$ 
  to have a permutation of $\{ 1, \ldots, n + 1 \}$ that has $k$ fixed points. 
  \item If $i \neq j$ then we can obtain a correspondence between $\sigma$ and a permutation $\sigma'$ of
  $\{ 1, \ldots, n \}$ from $\sigma$ as follows: $\sigma' ( j ) = i$ and
  $\sigma' ( t ) = \sigma ( t )$ for $t \neq j$. It is clear that $\sigma'$ has $k$ fixed
  points. Conversely, for each permutation $\sigma'$ of $\{ 1, \ldots,
  n \}$ that has $k$ fixed points, we can choose any $j$ such that $\sigma' (
  j ) = i \neq j$, and get back a permutation $\sigma$ of $\{ 1, \ldots, n + 1
  \}$ that has $k$ fixed points by letting $\sigma ( t ) = \sigma' ( t )$ for $t
  \neq j, n + 1, \sigma ( j ) = n + 1$ and $\sigma ( n + 1 ) = \sigma' ( j ) =
  i$. 
\end{enumerate}
\end{enumerate}
Hence $f_{n + 1} ( k ) = f_n ( k - 1 )  + ( n - k ) f_n ( k ) + ( k + 1 ) f_n ( k + 1 )$ for all $n, k$. This
  concludes the proof.         
\end{proof}

Lemma \ref{re} can be applied to obtain the following identity for $f_n(k)$ (Note that $f_n ( k ) = 0$ whenever $k < 0$ or $k >n$ so we do not need to specify the range of $k$).

\begin{lemma} \label{2} Suppose that $n, k, t$ are integers, $t\geq -1$. Then
  \begin{equation*} \sum_{k} [k]_{t+1} f_n ( k ) = 
	\begin{cases}
    n!  & \text{if}\; \; n \geqslant t+1,\\
    0  & \text{otherwise}
   \end{cases}
  \end{equation*}  
\end{lemma}

\begin{proof}
  We prove this using a double induction. The outer induction is on $t$ and the inner
  one is on $n$. By convention, $[k]_0=1£$. Also we have $\sum_k f_n ( k ) = n!$ which is trivial from the definition 
  of $f_n ( k )$. Hence the claim holds for $t=-1$. Next, suppose that the claim holds
  for $t - 1$. We prove that it holds for $t$. Define
   \[F ( n, t ) := \sum_k [k]_{t+1} f_n(k) = \sum_k k ( k -1 ) \ldots ( k - t ) f_n ( k ).\]
  
  Suppose that $n \leq t$. If $f_n ( k ) \neq 0$ then $0 \leqslant k
  \leqslant n \leqslant t$. But this implies that $k ( k - 1 ) \ldots ( k - t ) = 0$. Hence $F ( n, t ) = 0$  if $n \leqslant t$.  
  
  Suppose that $n = t + 1$. Then 
  \[F ( n, t ) = \sum_k k ( k -1 ) \ldots ( k - (n-1) ) f_n ( k ) = n! f_n(n) = n!\]
  since all but the last term of the sum equal zero. Hence the claim holds for $n = t+1$. For the inner induction, 
  suppose that $F ( n, t ) = n!$ for some $n \geq t + 1$. We will show that $F ( n + 1,
  t ) = ( n + 1 ) !$ using Lemma \ref{re}. We have  
  
\begin{align}
F &( n + 1, t ) = \sum k ( k - 1 ) \ldots ( k - t ) f_{n + 1} ( k )\nonumber\\
&= \sum k ( k - 1 ) \ldots ( k - t ) \left( f_n ( k - 1 ) + (
n - k ) f_n ( k ) + ( k + 1 ) f_n ( k_{} + 1 ) \right)\nonumber\\
&= n F ( n, t ) + \sum k ( k - 1 ) \ldots ( k - t ) \left( f_n ( k - 1 ) - k f_n ( k ) + ( k + 1 ) f_n ( k + 1 ) \right)\nonumber\\
&= n F ( n, t ) + \sum ( k + 1 ) k ( k - 1 ) \ldots ( k - t + 1 ) f_n ( k )\nonumber\\
& \,\,\,\,\,- \sum k ( k - 1 ) \ldots ( k - t ) k f_n ( k) + \sum ( k - 1 ) \ldots ( k - t - 1 ) k f_n ( k )\nonumber\\
&= n F ( n, t ) + \sum k ( k - 1 ) \ldots ( k - t + 1 )
[( k + 1 ) - ( k - t ) k + ( k - t ) ( k - t - 1 )] f_n ( k)\nonumber\\
&= n F ( n, t ) + \sum k ( k - 1 ) \ldots ( k - t + 1 )[ ( k + 1 ) - ( k - t ) ( t + 1 )] f_n ( k )\nonumber\\
&= n F ( n, t ) + \sum k ( k - 1 ) \ldots ( k - t + 1 )
[( t + 1 ) - ( k - t ) t] f_n ( k )\nonumber\\
&= n F ( n, t ) + ( t + 1 ) F ( n, t - 1 ) - t F ( n, t )\nonumber\\
&= n n! + ( t + 1 ) n! - t n!\label{moi}\\ 
&= ( n + 1 ) !.\nonumber
\end{align}                      
  
To see (\ref{moi}), note that the claim is true for $t-1$ by the outer induction. So $F(n,t-1) = n!$. Also $F(n,t) = n!$ by the inner inductive hypothesis. Hence the claim holds for $n + 1$. Therefore, it holds for every $n, t$. This concludes the proof of the lemma.
\end{proof}

\section{Proof of Theorem 1}

If $l>n$ then both sides of (\ref{hard}) equal zero. Hence we may assume that $l \leq n$. To have a permutation with exactly $k$ fixed points, we can first choose $k$ fixed points in $\binom{n}{k}$ ways. Then for each set of $k$ fixed points, we have $d(n-k)$ ways to arrange the $n-k$ remaining numbers such that we have no more fixed points. Hence 

\begin{equation}\label{support} f_n ( k ) = \binom{n}{k}d(n-k).\end{equation}
Substituting (\ref{support}) and $t = l - 1$ into Lemma \ref{2}, we have

\begin{equation}\label{4} \sum_k k ( k - 1 ) \ldots ( k - l + 1 ) \binom{n}{k} d (n-k) = n!.\end{equation}              But we have
\begin{align} k ( k - 1 ) \ldots ( k - l + 1 ) \binom{n}{k} & = k ( k- 1 ) \ldots ( k - l + 1 )\frac{n \ldots ( n - k + 1 )}{k!} \nonumber \\
&= n \ldots ( n - l + 1 ) \frac{( n - l ) \ldots ( n - k + 1 )}{( k - l ) !} \nonumber\\
&= n \ldots ( n - l + 1 ) \binom{n-l}{k-l}.\label{6}
\end{align}

Substituting (\ref{6}) into (\ref{4}), we obtain (\ref{hard}). This concludes the proof of the theorem.

\section{Proof of Theorem 2}

Now, let $g(x) = a_m x^m +\ldots+ a_0$ be any polynomial with integer coefficients. From (\ref{stirling}), we can rewrite $g(x)$ as 
\begin{equation*} 
	g(x) = \sum_{i=0}^m \left\{a_i \sum_{j=0}^i S(i,j) [x]_j \right\}.
\end{equation*}
Hence
\begin{align}
\sum_k g(k) f_n( k ) &= \sum_k  \left\{ \sum_{i=0}^m \left( a_i \sum_{j=0}^i S(i,j) [k]_j \right)   \right\} f_n ( k )\nonumber\\
                     &= \sum_{i=0}^m \left\{a_i \sum_{j=0}^i S(i,j) \left( \sum_k [k]_j f_n(k) \right) \right\}\nonumber\\
                     &= \sum_{i=0}^m \left\{a_i \sum_{j=0}^i S(i,j) F(n,j-1). \right\}\label{7}
\end{align}

From Lemma \ref{2}, $F(n,j-1) = n!$ for all $0 \leq j \leq n$. Aslo, from (\ref{bell}) $B_i = \sum_{j=0}^i S(i,j)$. Thus, (\ref{7}) implies that
\begin{align}
\sum_k g(k) f_n( k ) &= \sum_{i=0}^m \left\{a_i \sum_{j=0}^i S(i,j) n! \right\}\nonumber\\
										 &= \left\{\sum_{i=0}^m a_i B_i\right\} n!.\label{10}
\end{align}

Substituting $f_n(k) = \binom{n}{k} d(n-k)$ into (\ref{10}), we obtain (\ref{second}). This concludes the proof of the theorem.

\section{An application}

In this section, we will apply Theorem 2 to prove the irreducibility of the standard representation of symmetric groups. Let $G = S_n$ be the symmetric group on $X = \{ 1, \ldots, n \}$. Let $\mathbb{C}$ denote the complex numbers. Let $GL(d)$ stand for the group of all $d\times d$ complex matrices that are invertible with respect to multiplication. 

\begin{definition}
A matrix representation of a group $G$ is a group homomorphism
\[\rho : G \rightarrow GL(d).\] Equivalently, to each $g \in G$ is assigned $\rho(g) \in GL(d)$ such that
\begin{enumerate}
	\item $\rho(1) = I$, the identity matrix,
	\item $\rho(gh) = \rho(g) \rho(h)$ for all $g, h \in G$.
\end{enumerate}
\end{definition}
The parameter $d$ is called the \textit{degree} or \textit{dimension} of the representation and is denoted by $\deg(\rho)$.
All groups have the trivial representation of degree 1 which is the one sending every $g \in G$ to the matrix (1). We denote the trivial representation by $1$. An important representation of the symmetric group $S_n$ is the permutation representation $\pi$, which is of degree $n$. If $\delta \in S_n$ then we let $\pi(\delta) = (r_{i,j})_{n\times n}$ where
\begin{equation*}
r_{i,j} = 
\begin{cases}
1 & \text{if} \;\; \delta(j) = i,\\
0 & \text{otherwise.}
\end{cases}
\end{equation*}

\begin{definition}
Let $G$ be a finite group and let $\rho$ be a matrix representation of $G$. Then the character of $\rho$ is
\begin{equation*}\chi_{\rho}(g) = \text{tr}\ \rho(g),\end{equation*}
where tr denotes the trace of a matrix. 
\end{definition}

It is clear from Definition 2 that if $\delta \in S_n$ then
\begin{align*}
\chi_1(\delta) &= 1,\\
\chi_{\pi}(\delta) &= \text{number of fixed points of} \ \delta.
\end{align*} 

\begin{definition}
Let $\chi$ and $\phi$ be characters of a finite group $G$. Then
\[\left\langle \chi, \phi \right\rangle = \frac{1}{|G|} \sum_{g\in G} \chi(g) \phi(g^{-1}).\]
\end{definition}

A matrix representation $\rho$ of a group is called irreducible if $\left\langle \chi_{\rho}, \chi_{\rho} \right\rangle = 1$. Maschke's Theorem (see \cite{al1, al2}) states that every representation of a finite group having positive dimension can be written as a direct sum of irreducible representations. The permutation representation $\pi$ can be written as a direct sum of the trivial representation $1$ and another representation $\sigma$. The representation $\sigma$ is called the \emph{standard representation} of $S_n$. We have $\chi_{\pi} = \chi_1 + \chi_{\sigma}$ since for any $\delta \in S_n$ then $\pi(\delta) = 1(\delta) \oplus \sigma (\delta)$. Thus, for all $\delta \in S_n$ then
\begin{equation*}
\chi_{\sigma}(\delta) = (\text{number of fixed points of} \ \delta) -1.
\end{equation*}

Now we want to prove that $\sigma$ is irreducible. In other words, we need to show $\left\langle \chi_{\sigma}, \chi_{\sigma} \right\rangle = 1$, which is equivalent to
\begin{equation}\label{rep}
\sum_{k = 0}^n ( k - 1 )^2 f_n ( k ) = n! 
\end{equation}
Identity (\ref{rep}) can be obtained easily from Theorem 2 as follows.
\begin{align*}
\sum_k ( k - 1 )^2 f_n ( k ) &= \sum ( k^2 - 2k + 1 ) f_n ( k )\\
                             &= (2 - 2 + 1)n! = n!
\end{align*}
since $B_0 = B_1 = 1$ and $B_2=2$. This implies the irreducibility of standard representation of symmetric groups.

\section{Acknowlegement}

I would like to thank Dr. Catherine Greenhill for careful reading of the manuscript and for suggesting valuable improvements.

\bibliography{mybib}

\end{document}